\documentclass[10pt]{article}
\usepackage{amsmath, amssymb, amsthm, mathrsfs}
\usepackage{enumerate}
\usepackage{colordvi}
\usepackage{tikz}
\usepackage{enumerate}
\usepackage{graphicx}
\usepackage{hyperref}
\newtheorem{claim}{Claim}[section]
\newtheorem{theorem}[claim]{Theorem}
\newtheorem{lemma}[claim]{Lemma}

\usepackage{color}

\definecolor{Myred}{cmyk}{0.0,1.0,1.0,0.00}
\definecolor{Mypurple}{rgb}{0.5,0.0,0.5}
\definecolor{Mygreen}{rgb}{0,0.4,0}

\begin{document}

\title{Magnetic Dirichlet Laplacian on a perturbed twisted tube}

\author
{
Diana Barseghyan$^{1}$, Ricardo Abreu Blaya$^{2}$, Juan Bory-Reyes$^{3}$,\\ and Baruch Schneider$^{1}$
}
\date{
\small $^{1}$ Department of Mathematics, University of Ostrava,  30.dubna 22, 70103 Ostrava, Czech Republic,
\\ 
\small $^{2}$ Facultad de Matem\'{a}ticas, Universidad Aut\'{o}noma de Guerrero, Mexico,
\\
\small $^{3}$ ESIME-Zacatenco, Instituto Polit\'ecnico Nacional, M\'exico,\\ CDMX. 07738. M\'exico.
\\
E-mails:\; diana.schneiderova@osu.cz,  baruch.schneider@osu.cz, rabreublaya@yahoo.es, juanboryreyes@yahoo.com
}
\maketitle

\begin{abstract}
It is well known that the spectrum of the Dirichlet Laplacian for a compact perturbation of a three-dimensional, periodically twisted tube is unstable with respect to domain deformations. This means that if the periodically twisted tube is unperturbed, then the spectrum of the Dirichlet Laplacian is purely essential. On the other hand, the perturbation of this domain produces eigenvalues below the essential spectrum. This paper considers the Dirichlet-Laplace operator with a magnetic field. We explicitly prove that the spectrum of the magnetic Laplacian is stable under small and local deformations of the domain.
\end{abstract}

\bigskip



\section{Introduction} \label{s: intro}

Mathematical physics is the study of the relationship between the geometry of a region and the spectral properties of the operators that describe its dynamics. Over the last 25 years, a lot of work has been done on Dirichlet Laplacians in infinitely extended tubular regions. One particularly interesting class of problems concerns the properties of twisted tubes.

The fact that twisting of a non-circular tube gives rise to an effective repulsive interaction was first noted in \cite{CB96}, a proper mathematical meaning was given to this fact through an appropriate Hardy-type inequality \cite{EKK08}. On the other hand, while a periodic twist of an infinite straight tube raises the threshold of the essential spectrum, it is reasonable to expect that the local slowdown will act as an attraction. This could thus give rise to a discrete spectrum of the corresponding Dirichlet Laplacian; in \cite{EK05}, it was demonstrated that it is indeed the case. The effect has been further investigated, in particular, the paper \cite{BRKS09} analyzed the asymptotic distribution of the eigenvalues in case when the "slowdown" perturbation is infinitely extended and has prescribed decay properties.

In this paper we address another aspect of the spectral problem. Specifically, we study the magnetic Dirichlet Laplacian with a compactly supported magnetic field defined on a tube with a local perturbation of the periodically twisted one. 
We demonstrate that a compactly supported magnetic field can destroy the discrete spectrum but save the essential spectrum. The similar effect of the magnetic on the spectral behaviour have been demonstrated in the following works \cite{BBS24}, \cite{BBS25} and \cite{BE21}.  

\section{Magnetic Dirichlet Laplacian} \label{Our model}
\setcounter{equation}{0}

Let $\omega$ be an open, bounded domain in $\mathbb{R}^2$ such that its boundary is of class $C^2$, and let $\theta:\, \mathbb{R} \to \mathbb{R}$ be a differentiable function. For a given $s\in\mathbb{R}$ and $t:=(t_2, t_3)\in\omega$ we define the mapping $\mathfrak{L}:\: \mathbb{R}\times\omega\to\mathbb{R}^3$ by
\begin{equation}
\label{eq1}
\mathfrak{L}(s, t)=(s, t_2\cos\theta(s)+ t_3\sin\theta(s), t_3\cos\theta(s)- t_2\sin\theta(s))\,.
\end{equation}
The image $\mathfrak{L}(\mathbb{R}\times\omega)$ is a region in $\mathbb{R}^3$ that we call a \emph{twisted tube}
if the following two conditions are satisfied:

\begin{enumerate}[(i)]
 \setlength{\itemsep}{-3pt}

\item the function $\theta$ is not constant,

\item the set $\omega$ is not rotationally symmetric with respect to the origin in $\mathbb{R}^2$;

\end{enumerate}

\noindent in particular, the tube is set to be periodically twisted if $\theta$ is a linear function.

As we have indicated we are interested in tubes with a local perturbation of the periodic twisting. Consequently, we will  consider the angular function $\theta$ such that
\begin{equation}
\label{eq2}
\dot{\theta}(s)=\beta(s)= \beta_0-\mu(s)\,,
\end{equation}
where $\beta_0$ is a positive constant\footnote{The positivity assumption is made just for the sake of simplicity since for $\beta_0<0$, one can pass to a unitarily equivalent operator using a mirror-image transformation $(s, t)\mapsto (-s,t)$.} and $\mu(\cdot)$ is a positive and bounded function with $\mathrm{supp}\,\mu\subset[-s_0,s_0]$ for some $s_0>0$. We will be concerned with tubes
\begin{equation}
\label{Tube}
\Omega_\beta:=\mathfrak{L}(\mathbb{R}\times\omega)
\end{equation}
defined by $\mathfrak{L}= \mathfrak{L}_\theta$ corresponding to $\theta(s)= \int_{-s_0}^s \beta(s)\, \mathrm{d}s$, the angle  is determined up to a constant. 

In this paper, we will focus on the magnetic Dirichlet Laplacian, which is constructed as follows. Let us consider the closed quadratic form
$$Q_{\beta, \mathcal{A}}[u]=\|(i\nabla + \mathcal{A})u\|_{L^2(\Omega_\beta)},\quad u \in \mathcal{H}_0^1(\Omega_\beta),$$
where the real-value function $\mathcal{A}$ is a vector potential.

According to the first representation theorem \cite{K66}, there is a unique, nonnegative self-adjoint operator associated with this form:
\begin{equation}\label{1}H_{\Omega_\beta}(\mathcal{A}) = (i\nabla + \mathcal{A})^2.\end{equation}  
Throughout this work, we will suppose that the magnetic vector potential, $\mathcal{A}$, is compactly supported. 

In order to formulate the result, we require the following operator on $L^2(\omega)$,
$$
h_{\beta_0}=-\partial_{t_2}^2- \partial_{t_3}^2- \beta_0^2(t_2\partial_{t_3}-t_3\partial_{t_2})^2
$$
which due to the smoothness of the boundary of $\omega$ has the domain
$$
D(h_{\beta_0}) = \mathcal{H}_0^1(\omega)\cap \mathcal{H}^2(\omega)\,.
$$

Since the Dirichlet Laplacian $-\partial_{t_2}^2- \partial_{t_3}^2$, has a purely discrete spectrum on the bounded domain and the second term is positive, the spectrum of $h_{\beta_0}$ is purely discrete as well. We denote
\begin{equation}
\label{eqE}
E =\inf\sigma(h_{\beta_0})\,.
\end{equation}
As shown in \cite{K08}, that $E$ is greater than the ground state eigenvalue of the Dirichlet Laplacian on $\omega$, 
whenever $\omega$ is not rotationally symmetric. 

Then the following results for the Dirichlet Laplacian $H_{\Omega_\beta}(0)$ with zero magnetic field take place  \cite{EK05}:

\begin{theorem}
The spectrum of $H_{\Omega_\beta}(0)$ is purely absolutely continuous and covers the half-line $[E,\infty)$.
\end{theorem}

\begin{theorem} 
Assume that $\omega$ is not rotationally symmetric and that
\begin{equation} \label{ass}
\int_{-s_0}^{s_0}\, (\dot\theta^2(s)-\beta_0^2)\,d s <0\, ,
\end{equation}
where $\dot\theta(\cdot)$ is given by (\ref{eq2}). Then, the operator $H_{\Omega_\beta}(0)$ has at least one eigenvalue of finite
multiplicity below the threshold of the essential spectrum.
\end{theorem}
 
Below, we present the results showing the influence of the magnetic field on the spectral behaviour of the operator $H_{\Omega_\beta}(\mathcal{A})$.

\begin{theorem}\label{main1}
Let $\Omega_\beta$ be given in (\ref{Tube}) with $\beta$ described by (\ref{eq2}). Let $B= \mathrm{rot}(\mathcal{A}) \in C^1(\mathbb{R}^2)$ be a real-valued, compactly supported magnetic field. Then, the essential spectrum of the operator $(\ref{1})$ coincides with the half-line $[E, \infty)$.
\end{theorem}
 
\begin{theorem}\label{main}
Suppose the assumptions of Theorem \ref{main1}. Furthermore, suppose that $\omega$ contains the origin and its boundary is of class $C^2$. Let the magnetic vector potential be of the form $\mathcal{A}=(0, A, 0)$ where $A$ satisfies the following condition: there exists a two-dimensional ball, $\mathcal{B}(0, \tau)$, with center at zero and radius $\tau>0$ belonging to each cross-section of 
$\Omega_\beta\cap\{s= s^*\},\,s^*\in (-s_0,  s_0)$, such that the partial derivative $\partial_yA$ on $\mathcal{B}(0, \tau)$ does not depend on $s^*$ and is non-trivial as function of variable $t$. Then there exists a positive number, $\beta_0$, such that for
$\beta\in (0, 2\beta_0)$ the discrete spectrum of the operator $(\ref{1})$ below $E$ is empty.
\end{theorem}

\section{Proofs}

In this section, we will present the proofs of the above theorems. First, we will prove Theorem \ref{main}. Then, we will prove Theorem \ref{main1}.

\subsection{Proof of Theorem \ref{main}}

Given $\psi\in \mathcal{H}^1(\mathbb{R}\times\omega)$ it is useful to introduce the following shorthand,
$$\partial_\alpha \psi =t_2\partial_{t_3}\psi- t_3\partial_{t_2}\psi.$$

Then using the change of variables (\ref{eq1}) it is possible to check that for any $u\in \mathcal{H}_0^1(\Omega_\beta)$ we have

\begin{eqnarray*}
\partial_x u=\partial_{t_2}v\cos\theta+ \partial_{t_3}v\sin\theta,
\\
\partial_y u=-\partial_{t_2}v\sin\theta+ \partial_{t_3}v\cos\theta,
\\
\partial_s u=\partial_s v+ \dot{\theta}\partial_\alpha v\,,
\end{eqnarray*}
where 
\noindent \begin{equation}\label{representv} v(s, t)= u(s, t_2\cos \theta+ t_3\sin \theta, t_3\cos\theta- t_2\sin\theta)\end{equation}
\noindent and dot marks the derivative with respect to $s$. 

Hence with the notation \begin{equation}\label{notation}\tilde{A}(s, t)= A(s, t_2\cos\theta+ t_3\sin\theta, t_3\cos\theta- t_2\sin\theta)\end{equation} the quadratic form $Q_{\beta, \mathcal{A}}[u],\,u\in \mathcal{H}_0^1(\Omega_\beta)$ corresponding to operator (\ref{1}), can be rewritten as follows

\begin{eqnarray}\nonumber
Q_{\beta, \mathcal{A}}[u]- E\|u\|_{L^2(\Omega_\beta)}^2=\\ \nonumber 
\int_{\mathbb{R}\times\omega}\biggl(\left|i\partial _{t_2}v\cos\theta+i\partial _{t_3}v\sin\theta+\tilde{A}v\right|^2 +\left|-\partial_{t_2}v \sin\theta+\partial_{t_3}v\cos\theta\right|^2
+\left|\partial_sv+\dot{\theta}\partial _\alpha v\right|^2\biggr)\,d s\,d t\\ \nonumber
- E\|v\|_{L^2(\mathbb{R}\times\omega)}^2\\\nonumber
\end{eqnarray}

\begin{eqnarray}\nonumber
&=& \int_{\mathbb{R}\times\omega}\biggl(\left|
\partial_{t_2}v\right|^2 +i\tilde{A}\partial_{t_2}v\overline{v}\cos\theta
+ \left|\partial_{t_3}v\right|^2+ i\tilde{A}\overline{v}\partial_{t_3}v\sin\theta-  i\tilde{A}v\partial_{t_2}\overline{v}\cos\theta- i\tilde{A}v\partial_{t_3}\overline{v}\sin\theta \\ \nonumber 
&+&\tilde{A}^2 
|v|^2+ \left|\partial_sv\right|^2+ \dot{\theta}\partial_s v\partial_\alpha\overline{v}+\dot{\theta}\partial_s\overline{v} \partial_\alpha v+ \dot{\theta}^2\left|\partial_\alpha v\right|^2\biggr)\,d s\,d t -E\|v\|_{L^2(\mathbb{R}\times\omega)}^2\\\nonumber
&=& \int_{\mathbb{R}\times\omega}\left(\left|i\partial_{t_2}v+ \tilde{A}v\cos\theta \right|^2 +
\left|i\partial_{t_3}v+\tilde{A}v\sin\theta \right|^2\right)\,d s\,d t\\\label{straightform}
&+&\int_{\mathbb{R}\times\omega}\left( \left|\partial_s v\right|^2+ \dot{\theta}\partial_s v\partial _\alpha\overline{v}+ \dot{\theta}\partial_s\overline{v}\partial_\alpha v+ \dot{\theta}^2 \left|\partial_\alpha v\right|^2\right)\,d s\,d t
-E \|v\|_{L^2(\mathbb{R}\times\omega)}^2\,.
\end{eqnarray}

Using the fact that the ground-state eigenfunction $f$ of $h_{\beta_0}$ is strictly positive in $\omega$ --- cf.~\cite{EK05} --- we can decompose any $v\in C_0^\infty(\mathbb{R}\times\omega)$ as
\begin{equation}
\label{representvv}v(s, t)=f(t) g(s, t)\,.
\end{equation}
We substitute this into (\ref{straightform}) 
\begin{eqnarray}\nonumber 
&&Q_{\beta, \mathcal{A}}[u]- E\|u\|^2\\ \nonumber
&=&\int_{\mathbb{R}\times\omega}\biggl(|ig \partial_{t_2}f+if\partial_{t_2}g+
\tilde{A}fg\cos\theta|^2\biggr)\,d s\,d t
+\int_{\mathbb{R}\times\omega}\biggl(|ig\partial_{t_3}f+if\partial_{t_3}g+\tilde{A} fg\sin\theta|^2\biggr)\,d s\,d t
\\ \nonumber
&+&\int_{\mathbb{R}\times\omega}\biggl(f^2 |\partial_sg|^2+ \dot{\theta}f
\partial_sg(f \partial_\alpha\overline{g}+\overline{g}\partial_\alpha f)\biggr)\,d s\,d t+\int_{\mathbb{R}\times\omega}\dot{\theta}f\partial_s\overline{g}(f\partial_\alpha g+g\partial_\alpha f)\biggr)\,d s\,d t\\
\label{first expression} &+&
\int_{\mathbb{R}\times\omega}\dot{\theta}^2 |f\partial_\alpha g+g\partial_\alpha f|^2\,d s\,d t - E \int_{\mathbb{R}\times\omega}f^2|g|^2\,d s\,d t.
\end{eqnarray}

We will estimate the first five integrals on the right-hand side of (\ref{first expression}). We will start from the sum of first and second integrals as follows
\begin{eqnarray}\nonumber
&&\int_{\mathbb{R}\times\omega}\biggl(|ig\partial_{t_2}f+if\partial_{t_2}g+
\tilde{A}fg\cos\theta|^2\biggr)\,d s\,d t
+\int_{\mathbb{R}\times\omega}\biggl(|ig\partial_{t_3}f+if\partial_{t_3}g+\tilde{A}fg\sin\theta|^2\biggr)\,d s\,d t\\\nonumber &=&
\int_{\mathbb{R}\times\omega}\left(ig\partial_{t_2}f+if\partial_{t_2}g+\tilde{A}fg\cos\theta\right)
\left(-i\overline{g}\partial_{t_2}f-if\partial_{t_2}\overline{g}+\tilde{A}f\overline{g}\cos\theta\right)\\ \nonumber
&+&\int_{\mathbb{R}\times\omega}\left(ig \partial_{t_3}f+if\partial_{t_3}g+\tilde{A}fg\sin\theta\right)
\left(-i\overline{g} \partial_{t_3}f-if\partial_{t_3}\overline{g}+\tilde{A}f\overline{g}\sin\theta\right)
\\ \nonumber
&=&\int_{\mathbb{R}\times\omega}\biggl((\partial_{t_2}f)^2|g|^2+ fg\partial_{t_2}f \partial_{t_2}\overline{g}+f\overline{g}\partial_{t_2}g\partial_{t_2}f+|\partial_{t_2}g|^2f^2+ i\tilde{A}f^2 \overline{g}\partial_{t_2}g\cos\theta\\\nonumber
&-&i\tilde{A}f^2g\partial_{t_2}\overline{g}\cos\theta+\tilde{A}^2 f^2|g|^2\cos^2\theta+|g|^2(\partial_{t_3}f)^2+fg\partial_{t_3}\overline{g}\partial_{t_3}f+f\overline{g}\partial_{t_3}g \partial_{t_3}f\\ \label{*}
&+&f^2|\partial_{t_3}g|^2+i\tilde{A}f^2\overline{g}\partial_{t_3}g\sin\theta-i\tilde{A}f^2g \partial_{t_3}\overline{g}\sin\theta+\tilde{A}^2f^2|g|^2\sin^2\theta\biggr)\,d s\,d t\,.
\end{eqnarray}

Using integration by parts, it is easy to verify that
\begin{eqnarray*}
\int_{\mathbb{R}\times\omega}\left(fg\partial_{t_2}\overline{g}\partial_{t_2}f+f\overline{g}\partial_{t_2}g \partial_{t_2}f\right)\,d s\,d t
=-\int_{\mathbb{R}\times\omega}|g|^2(f\partial^2_{t_2}f+(\partial_{t_2}f)^2)\,d s\,d t\\
\int_{\mathbb{R}\times\omega}\left(fg\partial_{t_3}\overline{g}\partial_{t_3}f+f\overline{g}\partial_{t_3}g \partial_{t_3}f\right)\,d s\,d t
=-\int_{\mathbb{R}\times\omega}|g|^2 (f\partial^2_{t_3}f+(\partial_{t_3}f)^2)\,d s\,d t\,.\end{eqnarray*}

Using the above expressions in the right-hand side of (\ref{*}) one gets
\begin{eqnarray}\nonumber
&&\int_{\mathbb{R}\times\omega}\biggl(|ig\partial_{t_2}f+if\partial_{t_2}g+
\tilde{A}fg\cos\theta|^2\biggr)\,d s\,d t +\int_{\mathbb{R}\times\omega}\biggl(|ig\partial_{t_3}f+if\partial_{t_3}g+\tilde{A}fg\sin\theta|^2\biggr)\,d s\,d t\\\nonumber
&=&-\int_{\mathbb{R}\times\omega}|g|^2f(\partial_{t_2}^2f+\partial_{t_3}^2f)+f^2|\partial_{t_2}g|^2+i\tilde{A}f^2 \overline{g}\partial_{t_2}g\cos\theta-i\tilde{A}\partial_{t_2}\overline{g}gf^2\cos\theta\\ \nonumber
&+&\tilde{A}^2 f^2|g|^2\cos^2\theta
+|\partial_{t_3}g|^2 f^2+i\tilde{A}\overline{g}\partial_{t_2}gf^2\sin\theta-i\tilde{A}f^2\partial_{t_3}\overline{g} g \sin\theta
+\tilde{A}^2f^2|g|^2\sin^2\theta\biggr)\,d s\,d t\\ \label{magn.form}
&=&-\int_{\mathbb{R}\times\omega}|g|^2f(\partial_{t_2}^2f+\partial_{t_3}^2f)\,d s\,d t\\ \nonumber 
&+&\int_{\mathbb{R}\times\omega}
f^2 \left(|i\partial_{t_2}g+\tilde{A}g\cos\theta|^2+|i\partial_{t_3}g+\tilde{A}g\sin\theta|^2\right)\,d s\,d t\,.
\end{eqnarray}

Next, we will study the third and fourth integrals. Since
\begin{eqnarray*}
\int_{\mathbb{R}\times\omega}\left(f\overline{g}\partial_sg\partial_\alpha f+fg\partial_s\overline{g}\partial_\alpha f\right)\,d s\,d t=0\,, \\
\int_{\mathbb{R}\times\omega}\left(f\overline{g}\partial_\alpha g\partial_\alpha f+fg\partial_\alpha f\partial_\alpha \overline{g}\right)\,d s\,d t=-\int_{\mathbb{R}\times\omega}|g|^2\left(f \partial_\alpha^2f+(\partial_\alpha f)^2\right)
\,d s\,d t\,,
\end{eqnarray*}
\noindent  then using (\ref{eq2}) it is easy to check that
\begin{eqnarray}\nonumber
&&\int_{\mathbb{R}\times\omega}\biggl(f^2|\partial_sg|^2+ \dot{\theta}f
\partial_sg (f\partial_\alpha\overline{g}+\overline{g}\partial_\alpha f)\biggr)\,d s\,d t+\int_{\mathbb{R}\times\omega}\dot{\theta}f\partial_s\overline{g}(f\partial_\alpha g+g\partial_\alpha f)\biggr)\,d s\,d t\\ \nonumber
&+&\int_{\mathbb{R}\times\omega}\dot{\theta}^2|f\partial_\alpha g+g\partial_\alpha f|^2\,d s\,d t\\ \nonumber
&=& \int_{\mathbb{R}\times\omega}\biggl(f^2|\partial_sg|^2+ (\beta- \mu)f
\partial_sg (f\partial_\alpha\overline{g}+\overline{g}\partial_\alpha f)\biggr)\,d s\,d t+\int_{\mathbb{R}\times\omega}(\beta_0- \mu)f\partial_s\overline{g}(f\partial_\alpha g+g\partial_\alpha f)\biggr)\,d s\,d t\\ \nonumber
&+&\int_{\mathbb{R}\times\omega}(\beta_0- \mu)^2(f^2|\partial_\alpha g|^2+ f\overline{g}\partial_\alpha g \partial_\alpha f+  f g\partial_\alpha \overline{g}\partial_\alpha f+ |g|^2
(\partial_\alpha f)^2 )\,d s\,d t
\\ \label{last}
&=&\int_{\mathbb{R}\times\omega}f^2 |\partial_sg+\beta_0\partial_\alpha g|^2\,d s\,d t- \beta_0^2\int_{\mathbb{R}\times\omega}f|g|^2\partial_\alpha^2f- I_\mu\,,
\end{eqnarray}
where
\begin{eqnarray*}
I_\mu:=\int_{\mathbb{R}\times\omega}\biggl(\mu\,f^2\partial_sg\partial_\alpha\overline{g} +\mu\,f\overline{g} \partial_sg \partial_\alpha f+ \mu f^2\partial_s\overline{g}\partial_\alpha g+\mu\,fg\partial_s\overline{g} \partial_\alpha f\\+ (2\beta_0 \mu-\mu^2)(f^2|\partial_\alpha g|^2+ f\overline{g} \partial_\alpha g\partial_\alpha f+ fg \partial_\alpha f\partial_\alpha \overline{g}+ |g|^2(\partial_\alpha f)^2)\biggr)\,d s\,d t
\end{eqnarray*}

Combining (\ref{magn.form}) and (\ref{last}) and using the fact that $f$ is the ground state eigenfunction of operator $h_{\beta_0}$ corresponding to eigenvalue $E$, expression (\ref{first expression}) implies

\begin{eqnarray}\nonumber 
&&Q_{\beta, \mathcal{A}}[u]-E \|u\|^2\\ \nonumber
&=&-\int_{\mathbb{R}\times\omega}|g|^2 f(\partial_{t_2}^2f+\partial_{t_3}^2f)\,d s\,d t+\int_{\mathbb{R}\times\omega}
f^2 \left(|i\partial_{t_2}g+\tilde{A}g\cos\theta|^2+|i\partial_{t_3}g+\tilde{A}g\sin\theta|^2\right)\,d s\,d t\\ \nonumber
&+&\int_{\mathbb{R}\times\omega}f^2 |\partial_sg+ \beta_0\partial_\alpha g|^2\,d s\,d t- \beta_0^2 \int_{\mathbb{R}\times\omega}f|g|^2\partial_\alpha^2f- E\int_{\mathbb{R}\times\omega} f^2|g|^2\,d s\,d t- I_\mu\\\nonumber
&=&\int_{\mathbb{R}\times\omega}
f^2 \left(|i\partial_{t_2}g+\tilde{A}g\cos\theta|^2+|i\partial_{t_3}g+\tilde{A}g\sin\theta|^2\right)\,d s\,d t\\ \label{expr.}
&+&\int_{\mathbb{R}\times\omega}f^2 |\partial_sg+ \beta_0\partial_\alpha g|^2\,d s\,d t- I_\mu\,.
\end{eqnarray}

Let's estimate it. First, we start with $I_\mu$. Using the Cauchy inequality, we get
\begin{eqnarray}
\nonumber
I_\mu\le 2\left(\int_{\mathbb{R}\times\omega}\mu f^2|\partial_\alpha g|^2\,d s\,d t\right)^{1/2}
\left(\int_{\mathbb{R}\times\omega}\mu f^2|\partial_sg|^2\,d s\,d t\right)^{1/2}\\
\nonumber
+2\beta_0\int_{\mathbb{R}\times\omega}\mu f^2|\partial_\alpha g|^2\,d s\,d t
+\int_{\mathbb{R}\times\omega}\mu\left(f\overline{g}\partial_sg\partial_\alpha f+ fg\partial_s\overline{g} \partial_\alpha f\right)
\,d s\,d t\\\label{I}+
\int_{\mathbb{R}\times\omega}(2\beta_0\mu- \mu^2)\left(f\overline{g}\partial_\alpha g\partial_\alpha f+ fg \partial_\alpha\overline{g}\partial_\alpha f+ |g|^2(\partial_\alpha f)^2\right)\,d s\,d t
\end{eqnarray}

Using integration by parts and noting that the function $f$ is the ground state-eigenfunction of the operator $h_{\beta_0}$ we arrive
\begin{eqnarray*}
\int_{\mathbb{R}\times\omega}\mu\left(f\overline{g}\partial_sg\partial_\alpha f+ fg\partial_s\overline{g} \partial_\alpha f\right)
\,d s\,d t=-\int_{\mathbb{R}\times\omega}\dot{\mu}f\partial_\alpha f|g|^2\,d s\,d t\\= \frac{1}{2}\int_{\mathbb{R}\times\omega}\dot{\mu} 
(\overline{g}\partial_\alpha g+ g\partial_\alpha\overline{g}) f^2\,d s\,d t
\le \int_{\mathbb{R}\times\omega}|\dot{\mu}| 
f^2|g||\partial_\alpha g|\,d s\,d t\\\le \left(\int_{\mathbb{R}\times\omega} |\dot{\mu}|f^2|\partial_\alpha g|^2\,d s\,d t\right)^{1/2} \left(\int_{\mathbb{R}\times\omega} |\dot{\mu}|f^2|g|^2\,d s\,d t\right)^{1/2}\\\le \frac{1}{2}\int_{\mathbb{R}\times\omega} |\dot{\mu}|f^2|\partial_\alpha g|^2\,d s\,d t+ \frac{1}{2}\int_{\mathbb{R}\times\omega} |\dot{\mu}|f^2|g|^2\,d s\,d t\,,\end{eqnarray*}
\begin{eqnarray*}
&&\int_{\mathbb{R}\times\omega}(2\beta_0\mu- \mu^2)\left(f\overline{g}\partial_\alpha g\partial_\alpha f+ fg\partial_\alpha\overline{g}\partial_\alpha f+ |g|^2(\partial_\alpha f)^2\right)\,d s\,d t\\
&=&-\int_{\mathbb{R}\times\omega}( 2\beta_0 \mu- \mu^2)|g|^2f\partial_\alpha^2 f\,d s\,d t\\
&=&\frac{1}{\beta_0^2}\int_{\mathbb{R}\times\omega}( 2\beta_0 \mu- \mu^2)\left(E f+\partial^2_{t_2}f+\partial^2_{t_3}f\right) |g|^2 f\,d s\,d t
\end{eqnarray*}

\begin{eqnarray*}
&=& \frac{E}{\beta_0^2}\int_{\mathbb{R}\times\omega}( 2\beta_0\mu- \mu^2)f^2|g|^2\,d s\,d t+ \frac{1}{\beta_0^2}\sum_{j=2}^3 \int_{\mathbb{R}\times\omega}( 2\beta_0 \mu- \mu^2)\partial^2_{t_j} f |g|^2\,d s\,d t\\
&=&\frac{E}{\beta_0^2}\int_{\mathbb{R}\times\omega}( 2\beta_0\mu- \mu^2)f^2|g|^2\,d s\,d t- \frac{1}{\beta_0^2}\sum_{j=2}^3 \int_{\mathbb{R}\times\omega}( 2\beta_0 \mu- \mu^2)(\partial_{t_j} f)^2|g|^2\,d s\,d t\\
&-& \frac{1}{\beta_0^2}\sum_{j=2}^3 \int_{\mathbb{R}\times\omega}( 2\beta_0 \mu- \mu^2)f\partial_{t_j}f\left(\overline{g}\partial_{t_j}g+g\partial_{t_j}\overline{g}\right)\,d s\,d t\\ 
&\le&\frac{E}{\beta_0^2} \int_{\mathbb{R}\times\omega}( 2\beta_0 \mu- \mu^2)f^2|g|^2\,d s\,d t- \frac{1}{\beta_0^2}\int_{\mathbb{R}\times\omega}(2\beta_0 \mu- \mu^2)\left((\partial_{t_2}f)^2+(\partial_{t_3}f)^2\right)|g|^2\,d s\,d t\\
&+& \frac{2}{\beta_0^2} \sum_{j=2}^3\int_{\mathbb{R}\times\omega}(2\beta_0 \mu- \mu^2) f |g|\left|\partial_{t_j}f\right|\left|\partial_{t_j} g\right|\,d s\,d t\\
&\le& \frac{E}{\beta_0^2}\int_{\mathbb{R}\times\omega}(2\beta_0 \mu- \mu^2)f^2|g|^2\,d s\,d t-\frac{1}{\beta_0^2}\int_{\mathbb{R}\times\omega}(2\beta_0 \mu- \mu^2)\left((\partial_{t_2}f)^2+(\partial_{t_3}f)^2\right)|g|^2\,d s\,d t\\
&+& \frac{2}{\beta_0^2} \sum_{j=2}^3\left(\int_{\mathbb{R}\times\omega}(2\beta_0 \mu- \mu^2) f^2
|\partial_{t_j}g|^2\,d s\,d t\right)^{1/2}\left(\int_{\mathbb{R}\times\omega}(2\beta_0 \mu- \mu^2)|g|^2(\partial_{t_j}f)^2\,d s\,d t\right)^{1/2}
\\
&\le& \frac{E}{\beta_0^2}\int_{\mathbb{R}\times\omega}( 2\beta_0 \mu- \mu^2)f^2|g|^2\,d s\,d t-  \frac{1}{\beta_0^2} \int_{\mathbb{R}\times\omega}(2\beta_0 \mu- \mu^2)\left((\partial_{t_2}f)^2+(\partial_{t_3}f)^2\right) |g|^2\,d s\,d t\\
&+& \frac{1}{\beta_0^2}\int_{\mathbb{R}\times\omega}(2\beta_0 \mu- \mu^2)f^2\left(|\partial_{t_2}g|^2+ |\partial_{t_3}g|^2\right)\,d s\,d t+ \frac{1}{\beta_0^2}\int_{\mathbb{R}\times\omega}(2\beta_0 \mu- \mu^2)|g|^2\left((\partial_{t_2}f)^2+ (\partial_{t_3}f)^2\right)\,d s\,d t\\
&=& \frac{E}{\beta_0^2}\int_{\mathbb{R}\times\omega}( 2\beta_0 \mu- \mu^2)f^2|g|^2\,d s\,d t+ \frac{2}{\beta_0^2} \int_{\mathbb{R}\times\omega}(2\beta_0 \mu- \mu^2)f^2 \left(|\partial_{t_2}g|^2+|\partial_{t_3}g|^2\right)\,d s\,d t\,.\end{eqnarray*}

Applying the above expressions in the right-hand side of (\ref{I}) we have
\begin{eqnarray}\nonumber
I_\mu\le 2\left(\int_{\mathbb{R}\times\omega}\mu f^2|\partial_\alpha g|^2\,d s\,d t\right)^{1/2}
\left(\int_{\mathbb{R}\times\omega}\mu f^2|\partial_s g|^2\,d s\,d t\right)^{1/2}\\
\nonumber
+2\beta_0\int_{\mathbb{R}\times\omega}\mu f^2 |\partial_\alpha g|^2\,d s\,d t
+ \frac{1}{2}\int_{\mathbb{R}\times\omega} |\dot{\mu}|f^2|\partial_\alpha g|^2\,d s\,d t+ \frac{1}{2}\int_{\mathbb{R}\times\omega}|\dot{\mu}|f^2|g|^2\,d s\,d t\\\nonumber
+\frac{E}{\beta_0^2}\int_{\mathbb{R}\times\omega}( 2\beta_0 \mu- \mu^2)f^2|g|^2\,d s\,d t+ \frac{2}{\beta_0^2} \int_{\mathbb{R}\times\omega}(2\beta_0 \mu- \mu^2)f^2\left(|\partial_{t_2}g|^2+|\partial_{t_3}g|^2\right)\,d s\,d t
\\\nonumber 
\le c\int_{\mathbb{R}\times\omega}\mu f^2|\partial_s g|^2\,d s\,d t
+\int_{\mathbb{R}\times\omega} \left(c^{-1}\mu + 2\beta_0 \mu+\frac{1}{2}|\dot{\mu}| \right)\,f^2
|\partial_\alpha g|^2\,d s\,d t+ \frac{1}{2}\int_{\mathbb{R}\times\omega} |\dot{\mu}|f^2|g|^2\,d s\,d t
\\\nonumber
+\frac{E}{\beta_0^2}\int_{\mathbb{R}\times\omega}( 2\beta_0 \mu- \mu^2)f^2|g|^2\,d s\,d t+ \frac{2}{\beta_0^2} \int_{\mathbb{R}\times\omega}(2\beta_0 \mu- \mu^2)f^2\left(|\partial_{t_2}g|^2+|\partial_{t_3}g|^2\right)\,d s\,d t
\end{eqnarray}
for an arbitrary $c>0$ (to this number we will return a bit later). Further in view of the following obvious inequalities
\begin{eqnarray}\nonumber
|\partial_\alpha g|\le d\,\sqrt{|\partial_{t_2} g|^2+ |\partial_{t_3} g|^2}\,,
\\\label{pointwise}
|\partial_{t_2}g|^2+ |\partial_{t_3}g|^2\le 2(|i\partial_{t_2}g+\tilde{A}g\cos\theta|^2+|i\partial_{t_3}g+\tilde{A}g\sin\theta|^2)+ 2\tilde{A}^2|g|^2\,,
\end{eqnarray}
where 
$
d:=\sup_{(t_2,t_3)\in\omega}\sqrt{t_2^2+t_3^2},
$
and the fact that $\mathrm{supp}\,\mu\subset [-s_0, s_0]$ the above bound implies
\begin{eqnarray}\nonumber
&&I_\mu\le c\int_{(-s_0, s_0)\times\omega}\mu f^2|\partial_s g|^2\,d s\,d t\\ \nonumber
&+&2d^2\left(\frac{1}{c}+2\beta_0+\frac{1}{2}+\frac{4}{\beta_0 d^2}\right)\mathrm{max}\left\{\|\mu\|_\infty, \|\dot{\mu}\|_\infty\right\}\, \int_{(-s_0, s_0)\times\omega}f^2(|i\partial_{t_2}g+\tilde{A}g\cos\theta|^2\\ \nonumber
&+&|i\partial_{t_3}g+\tilde{A}g\sin\theta|^2)\,d s\,d t\\ \label{Imu}
&+&\left(2d^2\left(\frac{1}{c}+2\beta_0+\frac{1}{2}\right)\|\tilde{A}\|_\infty^2\right.\\ \nonumber 
&+&\left.\frac{2E}{\beta_0}+\frac{8}{\beta_0}\|\tilde{A}\|_\infty^2+\frac{1}{2}\right)
\int_{(-s_0, s_0)\times\omega}\mathrm{max}\left\{|\mu|, |\dot{\mu}|\right\} f^2|g|^2\,d s\,d t\,.
\end{eqnarray}

Let us now return to the quadratic form expression (\ref{expr.}). If
$$
2\beta_0\int_{(-s_0, s_0)\times\omega}f^2|\partial_sg||\partial_\alpha g|\,d s\,d t\le\frac{1}{2}
\int_{(-s_0, s_0)\times\omega}f^2|\partial_sg|^2\,d s\,d t
$$
then it is easy to check that estimating the first two integrals in the right-hand side of expression (\ref{expr.}) one obtains
\begin{eqnarray}
\nonumber Q_{\beta, \mathcal{A}}[u]- E \|u\|^2\\\nonumber\ge\frac{1}{2}\int_{(-s_0, s_0)\times\omega}
f^2 \left(|i\partial_{t_2}g+\tilde{A}g\cos\theta|^2+|i\partial_{t_3}g+\tilde{A}g\sin\theta|^2+|\partial_sg|^2\right)\,d s\,d t- I_\mu\,.
\end{eqnarray}

Hence the said formula and (\ref{Imu}) yield the estimate
\begin{eqnarray}\nonumber
Q_{\beta, \mathcal{A}}[u]- E\|u\|^2\ge\\\nonumber \gamma_{\beta_0}^1\int_{(-s_0, s_0)\times\omega} f^2 \left(|i \partial_{t_2}g+\tilde{A} g\cos\theta|^2+|i\partial_{t_3}g+\tilde{A}g\sin\theta|^2+|\partial_s g|^2\right)\,d s\,d t\\\label{Qpsi}-\gamma_{\beta_0}^2\int_{(-s_0, s_0)\times\omega}\mathrm{max}\left\{|\mu|, |\dot{\mu}|\right\} f^2|g|^2\,d s\,d t,
\end{eqnarray}
where
$$
\gamma^1_{\beta_0}= \mathrm{max}\left\{\frac{1}{2}-2d^2\left(\frac{1}{c}+2\beta_0+\frac{1}{2}+\frac{4}{\beta_0 d^2}\right)\mathrm{max}\left\{\|\mu\|_\infty, \|\dot{\mu}\|_\infty\right\}, \frac{1}{2}-c\|\mu\|_\infty\right\}\,,
$$
\begin{equation}\label{gamma2}
\gamma^2_{\beta_0}=2d^2\left(\frac{1}{c}+2\beta_0+\frac{1}{2}+\frac{4}{\beta_0 d^2}\right)\|\tilde{A}\|_\infty^2+\frac{2E}{\beta_0}+\frac{1}{2}\,.
\end{equation}

In the opposite case, namely under the assumption
$$
2\beta_0 \int_{(-s_0, s_0)\times\omega}f^2|\partial_sg||\partial_\alpha g|\,d s\,d t >\frac{1}{2}
\int_{(-s_0, s_0)\times\omega}f^2|\partial_sg|^2\,d s\,d t\,,
$$
by virtue of inequalities (\ref{pointwise}) and the Cauchy inequality applied to the left-hand side of the above bound we obtain  
\begin{eqnarray*}
\int_{(-s_0, s_0)\times\omega}f^2|\partial_sg|^2\,d s\,d t<32\beta_0^2d^2\int_{(-s_0, s_0)\times\omega}f^2
\left(|i\partial_{t_2}g+\tilde{A}g\cos\theta|^2+|i\partial_{t_3}g+\tilde{A}g\sin\theta|^2\right)\,d s\,d t\\+ 32\beta_0^2d^2\int_{(-s_0, s_0)\times\omega}\tilde{A}^2f^2|g|^2\,d s\,d t\,.
\end{eqnarray*}
	
Combining this estimate and (\ref{Imu}), we infer from (\ref{expr.}) that
\begin{eqnarray}\nonumber
&&Q_{\beta, \mathcal{A}}[u]- E\|u\|^2\\ \nonumber
&\ge&\tilde{\gamma}^1_{\beta_0}\int_{(-s_0, s_0)\times\omega}f^2\left(|i\partial_{t_2}g+\tilde{A}g\cos\theta|^2+|i\partial_{t_3}g+\tilde{A}g\sin\theta|^2\right)\,d s\,d t\\ \label{Q}
&-& \gamma^2_{\beta_0}\int_{(-s_0, s_0)\times\omega}\mathrm{max}\left\{|\mu|, |\dot{\mu}|\right\}f^2|g|^2\,d s\,d t\\\nonumber &-& 32cd^2\beta_0^2\int_{(-s_0, s_0)\times\omega}\tilde{A}^2f^2|g|^2\,d s\,d t\,,
\end{eqnarray}
where
$$
\tilde{\gamma}^1_{\beta_0}= 1-32cd^2\beta_0^2-2d^2\left(\frac{1}{c}+2\beta_0+\frac{1}{2}+\frac{4}{\beta_0 d^2}\right)\mathrm{max}\left\{\|\mu\|_\infty, \|\dot{\mu}\|_\infty\right\}\,.
$$

Combining inequalities (\ref{Qpsi}) and (\ref{Q}) we have
\begin{eqnarray}\nonumber
&&Q_{\beta, \mathcal{A}}[u]- E\|u\|^2\\ \nonumber
&\ge&\overline{\gamma}^1_{\beta_0}\int_{(-s_0, s_0)\times\omega}f^2\left(|i\partial_{t_2}g+\tilde{A}g\cos\theta|^2+|i\partial_{t_3}g+\tilde{A}g\sin\theta|^2\right)\,d s\,d t\\ \label{Q1} 
&-& \gamma^2_{\beta_0}\int_{(-s_0, s_0)\times\omega}\mathrm{max}\left\{|\mu|, |\dot{\mu}|\right\}f^2|g|^2\,d s\,d t\\\nonumber
&-& 32cd^2\beta_0^2\int_{(-s_0, s_0)\times\omega}\tilde{A}^2f^2|g|^2\,d s\,d t\,,
\end{eqnarray}
where 
\begin{equation}\label{not.} \overline{\gamma}^1_{\beta_0}=\mathrm{min}\{\gamma^1_{\beta_0}, \tilde{\gamma}^1_{\beta_0}\}\,.\end{equation}

Let us fix a positive number $\delta$ (to be chosen later) and let 
$\omega_\delta=\left\{(t_2, t_3)\in\omega: \mathrm{dist}\left((t_2, t_3), \omega\right)\ge\delta\right\}.$
We are going to consider two different situations:

\begin{enumerate}[(a)]
\setlength{\itemsep}{3pt}
\item
$$\int_{[-s_0, s_0]\times\omega_\delta}f^2|g|^2\,d s\,d t\ge \frac{1}{2}\int_{[-s_0, s_0]\times\omega}f^2|g|^2\,d s\,d t\,,$$

\item
$$\int_{[-s_0, s_0]\times\omega_\delta}f^2|g|^2\,d s\,d t< \frac{1}{2}\int_{[-s_0, s_0]\times\omega}f^2|g|^2\,d s\,d t\,.$$
\end{enumerate}

Let the situation (a) takes place. We are going to employ the influence of the magnetic field. Before continuing, we refer the reader to the proof of Lemma 1 in the work \cite{BE21} where the authors prove that actually the ground state eigenvalue of the magnetic Neumann Laplacian on bounded two dimensional domain with magnetic field which has non-zero integral over some ball is lower bounded by strictly positive constant independently on the values of the magnetic field outside of the mentioned ball.

In view of (\ref{notation}) one has
\begin{eqnarray*}
\partial_{t_2}\tilde{A}= \partial_xA\cos\theta- \partial_yA\sin\theta\,,\\
\partial_{t_3}\tilde{A}= \partial_xA\sin\theta+ \partial_yA\cos\theta\,.
\end{eqnarray*}

Hence for any fixed $s\in(-s_0, s_0)$ the magnetic field corresponding to two dimensional vector potential $(\tilde{A}\cos\theta, \tilde{A}\sin\theta)$ defined on $\omega$ coincides with  $-\partial_yA$ defined on  
$\{(t_2\cos\theta(s)+ t_3\sin\theta(s), t_3\cos\theta(s)- t_2\sin\theta(s)): \,(t_2, t_2)\in \omega\}$.
Moreover, this quantity does not depend on $s\in(-s_0, s_0)$ and is non-trivial on ball $\mathcal{B}(0, \tau)$ due to assumptions of Theorem (\ref{main}). Hence in view of our above discussions 
there exists a constant $\tau_{A}>0$ such that the ground state eigenvalue of the Neumann magnetic Laplacian with vector potential $(\tilde{A}\cos\theta, \tilde{A}\sin\theta)$ defined on $(-s_0, s_0)\times\omega$ is greater than $\tau_A$.
Together with this fact and assumption $(a)$ (which means that \newline\noindent$\int_{(-s_0, s_0)\times\omega}f^2 |g|^2\,d s\,d t\le 2\int_{(-s_0, s_0)\times\omega_\delta}f^2 |g|^2\,d s\,d t$) inequality (\ref{Q1}) implies

\begin{eqnarray}\nonumber
Q_{\beta, \mathcal{A}}[u]- E\|u\|^2\\\nonumber\ge
\overline{\gamma}^1_{\beta_0}\tau_{A}\underset{\omega_\delta}{\mathrm{inf}}(f^2) \int_{(-s_0, s_0)\times\omega_\delta}|g|^2\,d s\,d t- \gamma^2_{\beta_0}\mathrm{max}\left\{\|\mu\|_\infty, \|\dot{\mu}\|_\infty\right\}\int_{(-s_0, s_0)\times\omega}f^2|g|^2\,d s\,d t\\\nonumber- 32cd^2\beta_0^2\int_{(-s_0, s_0)\times\omega}\tilde{A}^2f^2|g|^2\,d s\,d t\\\nonumber\ge\left(\frac{\overline{\gamma}^1_{\beta_0}\tau_{A}\underset{\omega_\delta}{\mathrm{inf}}(f^2)}{\|f\|_\infty^2}- 2\gamma^2_{\beta_0}\mathrm{max}\left\{\|\mu\|_\infty, \|\dot{\mu}\|_\infty\right\}\right)\int_{(-s_0, s_0)\times\omega_\delta}f^2|g|^2\,d s\,d t\\\nonumber- 64cd^2 \beta_0^2\|\tilde{A}\|_\infty^2\int_{(-s_0, s_0)\times\omega_\delta}f^2|g|^2\,d s\,d t
\,.\end{eqnarray}

Using the expressions (\ref{gamma2}) and (\ref{not.}) for the constants $\gamma^2_{\beta_0}$ and
$\overline{\gamma}^1_{\beta_0}$ and the fact of strictly positivity of the function $f$ \cite{EK05} one can easily check that by choosing number $c$ and function $\mu$ small enough it is possible to guarantee the validity of   
$$
\frac{\overline{\gamma}^1_{\beta_0}\tau_{A}\underset{\omega_\delta}{\mathrm{inf}}(f^2)}{\|f\|_\infty^2}- 2\gamma^2_{\beta_0}\mathrm{max}\left\{\|\mu\|_\infty, \|\dot{\mu}\|_\infty\right\}- 64cd^2\beta_0^2\|\tilde{A}\|_\infty^2> 0
$$ 
which means that
\begin{equation}
\label{1case}
Q_{\beta, \mathcal{A}}[u]- E\|u\|^2\ge0\,.
\end{equation}

Now we pass to case $(b)$. In view of inequality (\ref{Q1}) and the simple inequality 
$|a+ b|^2\ge \frac{|a|^2}{2}- |b|^2$, we have
\begin{eqnarray}\nonumber
&&Q_{\beta, \mathcal{A}}[u]- E\|u\|^2\\\nonumber
&\ge&\frac{\overline{\gamma}^1_{\beta_0}}{2}\int_{(-s_0, s_0)\times\omega}f^2\left(|\partial_{t_2}g|^2+|\partial_{t_3}g|^2\right)\,d s\,d t
-\overline{\gamma}^1_{\beta_0}\int_{(-s_0, s_0)\times\omega}\tilde{A}^2f^2 |g|^2\,d s\,d t\\ \label{q1}
&-&\gamma^2_{\beta_0}\int_{(-s_0, s_0)\times\omega}\mathrm{max}\left\{|\mu|, |\dot{\mu}|\right\}f^2|g|^2\,d s\,d t\\\nonumber &-& 32cd^2\beta_0^2\int_{(-s_0, s_0)\times\omega}\tilde{A}^2f^2|g|^2\,d s\,d t\,.
\end{eqnarray}

Let us return to the notation (\ref{representvv}): $v(s, t)=f(t) g(s, t)$. Doing the integration by parts, it is easy to verify that 
\begin{eqnarray}\nonumber
\int_{[-s_0, s_0]\times\omega}\left(|\partial_{t_2}v|^2+|\partial_{t_3}v)|^2\right)\,d s\,d t
\\\nonumber
=\int_{[-s_0, s_0]\times\omega}\left(|f \partial_{t_2}g+ \partial_{t_2}f g|^2\right)\,d s\,d t
+\int_{[-s_0, s_0]\times\omega}\left(|f \partial_{t_3}g+ \partial_{t_3}f g|^2\right)\,d s\,d t
\\\label{magn.}=\int_{[-s_0, s_0]\times\omega}f^2\left(|\partial_{t_2}g|^2+ |\partial_{t_3} g|^2\right)\,d s\,d t
-\int_{[-s_0, s_0]\times\omega}f \Delta f |g|^2\,d s\,d t\,.
\end{eqnarray}

To continue the proof, we will use the following facts:
\begin{enumerate}
\item[(i)]
According to standard elliptic regularity theory (see \cite[Sec.~6.3]{E98}), we have $f \in H^4(\omega)$.
Consequently, $\nabla f \in H^3(\omega)$ and $\Delta f \in H^2(\omega)$.
Using the Sobolev embedding \cite[Thm.~5.4]{A75}
$H^2(\omega) \hookrightarrow C^0(\overline{\omega})$.
Thus, we have that 
$\|\Delta f\|_{L^\infty(\omega)} < \infty$.
\item[(ii)]
For any domain~$\omega$ such that~$\partial\omega$ is of
class~$C^2$, there exists \cite[Lem.~4.6.1]{D89}
a positive number~$\alpha_0$ such that $f \geq \alpha_0 \tau$,
where $\tau(t):=\mathrm{dist}(t,\partial\omega)$.
\item[(iii)]
For any strongly regular domain~$\omega$, which is in particular
satisfied under the present smoothness assumption, the Hardy
inequality $-\Delta_D^\omega \geq c_0/\tau^2,\, c_0>0$, holds true
\cite[Sec.~1.5]{D89}.
\end{enumerate}
This, together with  (\ref{q1})  and (\ref{magn.}), gives
\begin{eqnarray}\nonumber 
Q_{\beta, \mathcal{A}}[u]- E\|u\|^2\ge\\ \nonumber\frac{\overline{\gamma}^1_{\beta_0}}{2}\left(\int_{[-s_0, s_0]\times\omega}\left(|\partial_{t_2}v|^2+|\partial_{t_3}v|^2\right)\,d s\,d t+ \int_{[-s_0, s_0]\times\omega}f \Delta f |g|^2\,d s\,d t\right)-\overline{\gamma}^1_{\beta_0}\int_{(-s_0, s_0)\times\omega}\tilde{A}^2 |v|^2\,d s\,d t\\\nonumber
\end{eqnarray}

\begin{eqnarray}\nonumber 
&-&\gamma^2_{\beta_0}\int_{(-s_0, s_0)\times\omega}\mathrm{max}\left\{|\mu|, |\dot{\mu}|\right\}|v|^2\,d s\,d t- 32cd^2\beta_0^2\int_{(-s_0, s_0)\times\omega}\tilde{A}^2|v|^2\,d s\,d t\\ \nonumber
&\ge&\frac{ \overline{\gamma}^1_{\beta_0}}{2}\left(\int_{[-s_0, s_0]\times\omega}\left(|\partial_{t_2}v|^2+|\partial_{t_3}v|^2\right)\,d s\,d t-  \|\Delta f\|_\infty\int_{[-s_0, s_0]\times\omega}\frac{|v|^2}{f}\,d s\,d t\right)\\ \nonumber
&-&\overline{\gamma}^1_{\beta_0}\int_{(-s_0, s_0)\times\omega}\tilde{A}^2
 |v|^2\,d s\,d t
-\gamma^2_{\beta_0}\int_{(-s_0, s_0)\times\omega}\mathrm{max}\left\{|\mu|, |\dot{\mu}|\right\}|v|^2\,d s\,d t\\\nonumber
&-& 32cd^2\beta_0^2\int_{(-s_0, s_0)\times\omega}\tilde{A}^2|v|^2\,d s\,d t\\ \nonumber
&\ge&\frac{ \overline{\gamma}^1_{\beta_0}}{2}\left(\int_{[-s_0, s_0]\times\omega}\left(|\partial_{t_2}v|^2+|\partial_{t_3}v|^2\right)\,d s\,d t-  \frac{\|\Delta f\|_\infty}{\alpha_0}\int_{[-s_0, s_0]\times\omega}\frac{|v|^2}{\tau(t)}\,d s\,d t\right)\\\nonumber
&-&\overline{\gamma}^1_{\beta_0}\int_{(-s_0, s_0)\times\omega}\tilde{A}^2
|v|^2\,d s\,d t
-\gamma^2_{\beta_0}\int_{(-s_0, s_0)\times\omega}\mathrm{max}\left\{|\mu|, |\dot{\mu}|\right\}|v|^2\,d s\,d t \\\nonumber
&-& 32cd^2\beta_0^2\int_{\mathbb{R}\times\omega}\tilde{A}^2|v|^2\,d s\,d t\\ \nonumber
&\ge&\frac{ \overline{\gamma}^1_{\beta_0}}{2}\left(\int_{[-s_0, s_0]\times\omega}\left(|\partial_{t_2}v|^2+|\partial_{t_3}v|^2\right)\,d s\,d t- \frac{\|\Delta f\|_\infty}{2\alpha_0}\int_{[-s_0, s_0]\times\omega}\left( \frac{\varepsilon |v|^2}{\tau(t)^2}+ \frac{1}{\varepsilon}|v|^2\right)\,d s\,d t\right)\\ \nonumber
&-&\overline{\gamma}^1_{\beta_0}\int_{(-s_0, s_0)\times\omega}\tilde{A}^2
|v|^2\,d s\,d t
-\gamma^2_{\beta_0}\int_{(-s_0, s_0)\times\omega}\mathrm{max}\left\{|\mu|, |\dot{\mu}|\right\}|v|^2\,d s\,d t\\\nonumber
&-& 32cd^2\beta_0^2\int_{(-s_0, s_0)\times\omega}\tilde{A}^2|v|^2\,d s\,d t\\ \nonumber
&\ge&\frac{ \overline{\gamma}^1_{\beta_0}}{2}\left(1- \frac{\varepsilon \|\Delta f\|_\infty}{2 c_0 \alpha_0}\right)\int_{[-s_0, s_0]\times\omega}\left(|\partial_{t_2}v|^2+|\partial_{t_3}v|^2\right)\,d s\,d t
-\frac{ \overline{\gamma}^1_{\beta_0} \|\Delta f\|_\infty}{4 \varepsilon\alpha_0}\int_{[-s_0, s_0]\times\omega}|v|^2\,d s\,d t\\ \label{q2}
&-&\overline{\gamma}^1_{\beta_0}\int_{(-s_0, s_0)\times\omega}\tilde{A}^2
|v|^2\,d s\,d t
-\gamma^2_{\beta_0}\int_{(-s_0, s_0)\times\omega}\mathrm{max}\left\{|\mu|, |\dot{\mu}|\right\}|v|^2\,d s\,d t\\\nonumber
&-& 32cd^2\beta_0^2\int_{(-s_0, s_0)\times\omega}\tilde{A}^2|v|^2\,d s\,d t\,,
\end{eqnarray}
where $\varepsilon$ is a positive number such that $1- \frac{\varepsilon \|\Delta f\|_\infty}{2 c_0 \alpha_0}> 0$.

To continue the proof, we need the following lemma (the proof we will give in the last section).

\begin{lemma}\label{lemma}
Suppose assumption (b) is valid. Then, there exists a constant $C>0$ such that 
$$
\int_{[-s_0, s_0]\times \omega}\left(|\partial_{t_2}v|^2+|\partial_{t_3}v|^2\right)\,d s\,d t\ge \frac{C}{\delta^2}\int_{[-s_0, s_0]\times \omega}|v|^2\,d s\,d t\,.$$
\end{lemma}

Using the above lemma and the compactness of the supports of the magnetic vector potential and the perturbation 
$\mu$  we can estimate the right-hand side of inequality (\ref{q2}) as follows

\begin{eqnarray}\nonumber Q_{\beta, \mathcal{A}}[u]-E\|u\|^2
\\\nonumber\ge
\frac{ \overline{\gamma}^1_{\beta_0} C}{2 \delta^2}\left(1- \frac{\varepsilon \|\Delta f\|_\infty}{2 c_0 \alpha_0}\right)\int_{[-s_0, s_0]\times\omega}|v|^2\,d s\,d t
\\\nonumber-\frac{ \overline{\gamma}^1_{\beta_0} \|\Delta f\|_\infty}{4 \varepsilon\alpha_0}\int_{[-s_0, s_0]\times\omega}|v|^2\,d s\,d t-\overline{\gamma}^1_{\beta_0}\int_{(-s_0, s_0)\times\omega}\tilde{A}^2
|v|^2\,d s\,d t
\\\nonumber-\gamma^2_{\beta_0}\int_{(-s_0, s_0)\times\omega}\mathrm{max}\left\{|\mu|, |\dot{\mu}|\right\}|v|^2\,d s\,d t- 32cd^2\beta_0^2\int_{(-s_0, s_0)\times\omega}\tilde{A}^2|v|^2\,d s\,d t\\\nonumber
\ge\biggl(\frac{ \overline{\gamma}^1_{\beta_0} C}{2 \delta^2}\left(1- \frac{\varepsilon \|\Delta f\|_\infty}{2 c_0 \alpha_0}\right)-\frac{ \overline{\gamma}^1_{\beta_0} \|\Delta f\|_{L^\infty(\omega)}}{4 \varepsilon\alpha_0}-\overline{\gamma}^1_{\beta_0}\|\tilde{A}\|_\infty^2
\\\label{ass1}-\gamma^2_{\beta_0}\mathrm{max}\left\{\|\mu\|_\infty, \|\dot{\mu}\|_\infty\right\}- 32cd^2\beta_0^2 \|\tilde{A}\|_\infty^2\biggr)\int_{(-s_0, s_0)\times\omega}|v|^2\,d s\,d t\,.
\end{eqnarray}

Let us choose $\delta$ in such a way that 
\begin{eqnarray*}
\nonumber
\frac{ \overline{\gamma}^1_{\beta_0} C}{2 \delta^2}\left(1- \frac{\varepsilon \|\Delta f\|_\infty}{2 c_0 \alpha_0}\right)-\frac{ \overline{\gamma}^1_{\beta_0} \|\Delta f\|_{L^\infty(\omega)}}{4 \varepsilon\alpha_0}-\overline{\gamma}^1_{\beta_0}\|\tilde{A}\|_\infty^2
\\-\gamma^2_{\beta_0}\mathrm{max}\left\{\|\mu\|_\infty, \|\dot{\mu}\|_\infty\right\}- 32cd^2\beta_0^2 \|\tilde{A}\|_\infty^2\ge 0\,,
\end{eqnarray*}
which guarantees the non-negativeness of the right-hand side of estimate (\ref{ass1}). 
This together with (\ref{1case}) finishes the proof of the theorem.

\subsection{Proof of Theorem \ref{main1}}

In this section, we will follow the proof method presented in \cite{BRKS09}.
The magnetic Laplacian operator $H_{\Omega_\beta}(\mathcal{A})$  can be written as
$$
H_{\Omega_\beta}(\mathcal{A}) = H_{\Omega_\beta}(0) + W,
$$
where
$$
W: =  2 i A\partial_{t_2}+ i (\partial_{t_2}A)+ A^2\,
$$
is the operator defined  on domain of $H_{\Omega_\beta}(\mathcal{A})$, that is on $\mathcal{H}_0^1(\Omega_\beta)$. Evidently, the operator $W (H_{\Omega_\beta}(0))^{-1/2}$ is bounded
in $L^2(\Omega_\beta)$. Due to the Sobolev embedding theorems
and the fact that $A$ is compactly supported it is easy to see that the operator $W (H_{\Omega_\beta}(0))^{-1}$ is compact in $L^2(\Omega_\beta)$. Therefore, the operator
$(H_{\Omega_\beta}(0))^{-1/2} W (H_{\Omega_\beta}(0))^{-1}$ is compact.
Since the operator $(H_{\Omega_\beta}(\mathcal{A}))^{-1/2} (H_{\Omega_\beta}(0))^{1/2}$ and,
hence, $(H_{\Omega_\beta}(\mathcal{A}))^{-1} (H_{\Omega_\beta}(0))^{1/2}$ is bounded, the resolvent difference
$(H_{\Omega_\beta}(\mathcal{A}))^{-1} - (H_{\Omega_\beta}(0))^{-1}$ is compact. By \cite[Theorem 4, Section 1, Chapter 9]{BS87}, we have
$$
\sigma_{\mathrm{ess}}(H_{\Omega_\beta}(\mathcal{A})) = \sigma_{\mathrm{ess}}(H_{\Omega_\beta}(0))= [E, \infty).
$$ 

\section{Proof of Lemma \ref{lemma}}

Since 
$$
\int_{[-s_0, s_0]\times \omega}\left(|\partial_{t_2}v|^2+|\partial_{t_3}v|^2\right)\,d s\,d t\ge \int_{[-s_0, s_0]\times (\omega\setminus\omega_\delta)}\left(|\partial_{t_2}v|^2+|\partial_{t_3}v|^2\right)\,d s\,d t
$$
it is enough to estimate  $\int_{[-s_0, s_0]\times (\omega\setminus\omega_\delta)}\left(|\partial_{t_2}v|^2+|\partial_{t_3}v|^2\right)\,d s\,d t$.

It is easy to notice that the set $\omega\setminus\omega_\delta$ can be presented as follows
\begin{equation}
\label{cupomega}
\omega\setminus\omega_\delta= \left(\cup_{i=1}^N\omega^i_\delta\right)\cup \left(\cup_{j=1}^M\omega^j_\delta\right)\,,
\end{equation}
where $N, M$ are some natural numbers and 
\begin{eqnarray}\nonumber \omega^i_\delta=\{t_2\in (x_1^i, x_2^i),\,t_3\in(y_1^i(t_2), y_2^i(t_2))\}, \\\label{omega}
\omega^j_\delta=\{t_3\in (y_1^j, y_2^j),\,t_2\in(x_1^j(t_3), x_2^j(t_3))\},
\end{eqnarray}
with some $x_1^i< x_2^i,\, \, y_1^j< y_2^j$ and smooth functions $x_1^j(t_3)< x_2^j(t_3)$, $y_1^i(t_2)< y_2^i(t_2)$ such that $\|x_2^j(t_3)- x_1^j(t_3)\|_{L^\infty(y_1^j, y_2^j)}\le C' \delta$ and
$\|y_2^i(t_2)- y_1^i(t_2)\|_{L^\infty(x_1^i, x_2^i)}\le C' \delta$ where $C'> 0$ is some constant depending on 
$\omega$.

Let us first estimate $\int_{(-s_0, s_0)\times\omega^i_\delta}|v|^2\,d s\,d t$ from below. Due to the fact that for any $t_2\in (x_1^i, x_2^i)$ we have $v(s, t_2, y_2^i(t_2))= 0,\,s\in(-s_0, s_0)$, then for each
$t_3\in (y_1^i(t_2), y_2^i(t_2))$ takes place
\begin{eqnarray*}
|v(s, t_2, t_3)|= \left|\int_{t_3}^{y_2^i(t_2)}\partial_z v(s, t_2, z)\,d z\right|\le (y_2^i(t_2)- t_3)^{1/2}\left(\int_{t_3}^{y_2^i(t_2)}|\partial_z v(s, t_2, z)|\,d z\right)^{1/2}\\\le (C' \delta)^{1/2} \left(\int_{t_3}^{y_2^i(t_2)}|\partial_z v(s, t_2, z)|^2\,d z\right)^{1/2}\le (C' \delta)^{1/2} \left(\int_{y_1^i(t_2)}^{y_2^i(t_2)}|\partial_z v(s, t_2, z)|^2\,d z\right)^{1/2}.
\end{eqnarray*}

Hence
\begin{eqnarray*}
\int_{(-s_0, s_0)\times\omega^1_\delta}|v|^2\,d s\,d t\le C' \delta \int_{(-s_0, s_0)\times\omega^1_\delta}\int_{y_1^i(t_2)}^{y_2^i(t_2)}|\partial_z v(s, t_2, z)|^2\,d z\,d s\,d t\\\le C' \delta \int_{(-s_0, s_0)\times\omega^1_\delta}(y_2^i(t_2)- y_1^i(t_2))|\partial_z v(s, t_2, z)|^2\,d z\,d s\,d t\\\le C'^2 \delta^2 \int_{(-s_0, s_0)\times\omega^1_\delta}|\partial_z v(s, t_2, z)|^2\,d z\,d s\,d t\,.\end{eqnarray*}

The rest cases with $\omega_\delta^j$, will be considered in the same way. This establishes Lemma \ref{lemma} with the constant $C=\frac{1}{C'^2}$.

\subsection*{Contributions}
Authors have been discussing and working together on the manuscript, contributing equally to the content, presentation, and reviewing the manuscript.
\bigskip

\subsection*{Data Availability}
No datasets were generated or analysed during the current study.
\bigskip

\subsection*{Competing interests}
The authors declare no competing interests.


\begin{thebibliography}{10}

\bibitem{A75} R. A.~ Adams, Sobolev spaces, Academic Press, New York, 1975.

\bibitem{BBS24} J.~ Bory-Reyes, D.~ Barseghyan, B.~ Schneider, Magnetic Schr\"{o}dinger operator with the potential supported in a curved two-dimensional strip, Mediterr. J. Math. 21(3) (2024), 1--15.

\bibitem{BBS25} D.~Barseghyan, J.~ Bory-Reyes, B.~Schneider, B.  Three-dimensional magnetic Schr\"{o}dinger operator with the potential supported in a tube, Annals of Functional Analysis 16(4) (2025).

\bibitem{BE21} D.~ Barseghyan, P.~ Exner, Magnetic field influence on the discrete spectrum of locally deformed leaky wires. Rep. Math. Phys. 88 (2021), 47--57.

\bibitem{BRKS09}
Ph.~Briet, G.~Raikov, H.~Kova\v{r}\'{\i}k, E.~Soccorsi,
Eigenvalue asymptotics in a twisted waveguide,
Comm. PDE 34 (2009), 818--836.

\bibitem{BS87} M. S.~Birman, M.Z.~Solomjak, Spectral Theory of Self-Adjoint Operators in
Hilbert Space, Reidel Publishing Company, Dordrecht, 1987.

\bibitem{CB96} I. J.~Clark, A. J.~Bracken, Bound states in tubular quantum waveguides with torsion,
J. Phys. A: Math. Gen. 29 (1996), 4527--4535.

\bibitem{D89} E. B.~ Davies, Heat kernels and spectral theory, Cambridge University Press,
1989.

\bibitem{E98}  L. C. Evans, Partial differential equations, Graduate Studies in 
Mathematics, vol. 19, American Mathematical Society, Providence, RI, 1998 

\bibitem{EK05}
P.~Exner, H.~Kova\v{r}\'{\i}k,  Spectrum of the Schr\"{o}dinger operator in a perturbed periodically twisted tube,
Lett. Math. Phys. 73 (2005), 183--192.

\bibitem{EKK08} T.~Ekholm, H.~Kova\v{r}\'{\i}k, D.~Krej\v{c}i\v{r}\'{\i}k,
A Hardy inequality in twisted waveguides, Arch. Rational Mech. Anal. 188 (2008), 245--264.

\bibitem{K66}
T.~Kato,  Perturbation Theory for Linear Operators, Classics in Mathematics, Springer-Verlag, Berlin 1995.

\bibitem{K08}  D.~Krej\v{c}i\v{r}\'{\i}k, Twisting versus bending in quantum waveguides,
Analysis on Graphs and Applications (Cambridge 2007), Proc. Sympos. Pure Math., vol. 77, pp. 617-636, Amer. Math. Soc., Providence, RI, 2008.

\end{thebibliography}
\end{document}